\documentclass{commat}

\title{%
   On the convergence of random Fourier-Jacobi series of continuous functions
    }

\author{
    Partiswari Maharana and Sabita Sahoo
    }

\affiliation{
    \address{Partiswari Maharana --
    Department of Mathematics, Sambalpur University, Odisha, India.
        }
    \email{
    partiswarimath1@suniv.ac.in
    }
    \address{Sabita Sahoo --
   Department of Mathematics, Sambalpur University, Odisha, India.
        }
    \email{
     sabitamath@suniv.ac.in
    }
    }

\abstract{%
The interest in orthogonal polynomials and random Fourier series in numerous branches of science and a few studies on random Fourier series in orthogonal polynomials inspired us to focus on random Fourier series in Jacobi polynomials.
In the present note, an attempt has been made to investigate the stochastic convergence of some random Jacobi series. We looked into the random series $\sum_{n=0}^\infty d_n r_n(\omega)\varphi_n(y)$ in orthogonal polynomials $\varphi_n(y)$ with random variables $r_n(\omega).$ The random coefficients $r_n(\omega)$ are the Fourier-Jacobi coefficients of continuous stochastic processes such as symmetric stable process and Wiener process. The $\varphi_n(y)$ are chosen to be the Jacobi polynomials and their variants depending on the random variables associated with the kind of stochastic process.
The convergence of random series is established for different parameters $\gamma,\delta$ of the Jacobi polynomials with corresponding choice of the scalars $d_n$ which are Fourier-Jacobi coefficients of a suitable class of continuous functions.
The sum functions of the random Fourier-Jacobi series associated with continuous stochastic processes are observed to be the stochastic integrals. The continuity properties of the sum functions are also discussed.
}



\VOLUME{32}
\NUMBER{1}
\YEAR{2024}
\firstpage{49}
\DOI{https://doi.org/10.46298/cm.10412}

\begin{paper}
\section{Introduction}

In the early 1800s, the Fourier series was invented to solve the problem of heat diffusion in a continuous medium.
It has widespread application in different branches of science. The beginning of the study of random series of functions in the 1930s by Paley and Zygmund in their pioneering works gave a new direction to look into Fourier series with random coefficients.
Kahane, Marcus, and Pisier have synthesized the most prominent work on the random Fourier series, and many problems are still open today.
The random Fourier series plays a vital role in electrical engineering, signal processing, optics, etc., because white noise, which is an inherent part of all these, is a random signal. It is also used in image encryption and decryption \cite{LL}, \cite{LL1}, and the generation of random noise \cite{GA}.

The orthogonal polynomials, which originated in the 19th century,
play an essential role in mathematical physics.
There has been extensive study of the Fourier series in orthogonal polynomials.
Marian and Marian \cite{MM} studied power series involving orthogonal polynomials, which occur in the problems of quantum optics.
However, the random Fourier series in orthogonal polynomials have received very little attention.
Some of the reasons are the difficulty of dealing with it and the lack of literature on its application to physical problems.
In 2007, Liu and Liu  \cite{LL}, \cite{LL1} defined the random Fourier transform in Hermite polynomials and applied it in the field of image encryption and decryption. They also expected its other kind of applications in optics and information technology, etc. This motivated us to explore random Fourier series in orthogonal polynomials.

Liu and Liu obtained a method to define the random Fourier transform while investigating the multiplicity and complexity of eigenvalues of the fractional Fourier transform.
They used the Fourier series expansion in Hermite polynomials of functions $f$ that belong to $L^2(\mathbb{R})$ in the form of
\begin{equation*}
f(y):=\sum\limits_{n=1}^\infty d_n \varphi_n(y),
\end{equation*}
where
\begin{equation*}
d_n:=\int_{-\infty}^\infty f(y)\varphi_n(y) dy
\end{equation*}
 are the Fourier-Hermite coefficients of $f$ and  $\varphi_n(y)$ are the Hermite Gaussian functions.
Since Hermite Gaussian functions are eigenfunctions of the
 Fourier transform defined on $L^2(\mathbb{R}),$
with corresponding eigenvalues
\begin{equation*}
\lambda_n:=\exp\Big[\frac{-i \pi}{2}\mod(n, 4)\Big],
\end{equation*}
where$\mod(n, 4)$ is the set of integers modulo 4.
It leads to finding the Fourier transform $\mathfrak{F}$ of $f$ as
\begin{equation*}
\sum_{n=1}^\infty d_n \lambda_n \varphi_n(y).
\end{equation*}
Further, the eigenfunctions of fractional Fourier transform of $f$ of rational order $\beta$ are the same as the Hermite Gaussian functions, but correspond to the eigenvalues $\lambda_n^\beta.$
This helps to express the fractional Fourier transform $\mathfrak{F}^\beta$  of $f$ of rational order $\beta$ in series form as
\begin{equation}\label{0.3.1}
\sum_{n=0}^\infty d_n \lambda_n^\beta \varphi_n(y),
\end{equation}
where
\begin{equation}\label{0.7}
\lambda_n^\beta=: \exp\Big\{\frac{-i n\beta \pi}{2}\Big\},\; n=0,1,2,\dots.
\end{equation}
However $\lambda_n^\beta$ can also be expressed as
\begin{equation}\label{0.9}
\lambda_n^\beta:={\Big[\exp \displaystyle{\frac{-in \pi}{2}} \Big]}^\beta,\; n=0,1,\dots.
\end{equation}
which is same as the expression in (\ref{0.7}). But (\ref{0.9}) is different from (\ref{0.7}) as it generates many different values in comparison to the 4 fixed values obtained from (\ref{0.7}).
Suppose the rational number $\beta={\frac{P}{Q}},$ where $P$ and $Q$ are integers, then the fractional eigenvalue is written as
\begin{align}
\lambda_n^\beta
&:= {\Big[\exp{\frac{-i \pi n}{2}}\Big]}^\beta \nonumber \\
&= \exp{\Big[ -iP\Big(\frac{\pi n}{2} +2k \pi\Big)/Q \Big]} \label{0.8},
\end{align}
where $k=0,1,2,\dots, Q-1.$ Hence, there exist $Q$ possible eigenvalues $\lambda_n^\beta$ for every $n.$ For the $N$th partial sum of the series (\ref{0.3.1}), the eigenvalues $\lambda_n^\beta$ have $Q^N$ choices, which is a large number. Thus their are $Q^N$ different results for the $N$th partial
sum of the fractional Fourier transform. Now, if the fractional order $\beta$ is extended to an irrational order $r,$ let it be approximated by the fraction ${\frac{P}{Q}}.$ When the distance between $r$ and ${\frac{P}{Q}}$ approaches to zero, then the denominator $Q$ will approach to infinity and the number of eigenvalues in (\ref{0.8}) will spread all over the entire unit circle. Thus the eigenvalues corresponding to every eigenfunction can have infinite choices. The infinite number of eigenvalues make the choice of eigenvalues in (\ref{0.8}) in an absolutely random way from the unit circle in $\mathbb{C}.$ Let these infinite number of eigenvalues be denoted as
$ \mathcal{R}(\lambda_n),$ and the series $(\ref{0.3.1})$ become
\begin{equation}\label{0.5}
\mathcal{R}[f(y)]:=\sum_{n=0}^\infty d_n \mathcal{R}(\lambda_n) \varphi_n(y),
\end{equation}
where $\mathcal{R}(\lambda_n):=\exp[i\pi \; \mathrm{Random}(n)]$ are randomly chosen values on the unit circle in $\mathbb{C}.$
They called (\ref{0.5}) as the random Fourier transform instead of Fourier transform of irrational order.
However, this random Fourier transform is found to be a random Fourier series in orthogonal Hermite polynomials.
It raises the question of what would happen if the random coefficients $\mathcal{R}(\lambda_n)$ selected from the unit circle in $\mathbb{C}$ were replaced by other random variables in the random series (\ref{0.5}). 

Our purpose in this article is to look at the random series
\begin{equation}\label{0.4}
\sum_{n=0}^\infty d_n r_n(\omega)\varphi_n(y)
\end{equation}
in orthogonal polynomials $\varphi_n(y)$ with random coefficients $r_n(\omega)$ and the scalars $d_n.$ In this work, the orthogonal polynomials $\varphi_n(y)$ are considered to be the Jacobi polynomials.
The Jacobi
polynomials have great importance in a hierarchy of the orthogonal polynomial
classes and found widespread use in all areas of science and engineering.
Recently, Arenas, Ciaurri, and Labarga \cite{ACL} studied the convergence of the Fourier-Jacobi series in a discrete version.
For  parameters $\gamma,\delta> -1,$ denote $p_n^{(\gamma,\delta)}(y),\;n \in \mathbb{N}$ to be the orthonormal Jacobi polynomials with respect to the Jacobi weight
$\rho^{(\gamma,\delta)}(y):={(1-y)}^\gamma{(1+y)}^\delta,\;y \in[-1,1].$

The literature on random Fourier series (\ref{0.4}) gives a way to deal with random Fourier series in Jacobi polynomials.
Nayak, Pattanayak, Mishra \cite{NPM} discussed the convergence, summability of random Fourier series associated with symmetric stable process
of index $\alpha \in[1,2]$ and shown that the random Fourier series converges in probability to a stochastic integral.
We know that, if $X(t,\omega),\; t\in \mathbb{R}$ is a continuous stochastic process with independent increments and $f$ is a continuous function in $[a,b],$ then the stochastic integral
\begin{equation}\label{0.1}
\int\limits_a^bf(t)dX(t,\omega)
\end{equation}
is defined in the sense of probability and is a random variable (Lukacs \cite{L}).
 Further, if $X(t,\omega)$ is a symmetric stable process of index $\alpha \in[1,2],$ then the stochastic integral~(\ref{0.1}) is defined in the sense of probability, for $f \in L_{[a,b]}^p, \; p \geq \alpha \geq 1$ \cite{NPM}.
In particular, if $f(t)=\varphi_n(t)$ is the orthonormal Jacobi polynomial $p_n^{(\gamma,\delta)}(t),\;\gamma,\delta>-1,$ then $p_n^{(\gamma,\delta)}(t)\rho^{(\eta,\tau)}(t)$ is continuous in $[-1,1],$ for $\eta,\tau \geq 0$ and
hence the integrals
\begin{equation}\label{1.13.0}
A_n(\omega):=\int_{-1}^1 p_n^{(\gamma,\delta)}(t)\rho^{(\eta,\tau)}(t)dX(t,\omega)
\end{equation}
exist and are random variables, for either of the choice of the stochastic process $X(t,\omega).$ The random variables
$A_n(\omega)$ are called the Fourier-Jacobi coefficients of the corresponding stochastic process $X(t,\omega).$ These $A_n(\omega)$ can be seen to be independent
if $X(t,\omega)$ is the Wiener process but are no more independent if $X(t,\omega)$ is the symmetric stable process.
In this article the random coefficients $r_n(\omega)$ in the random series (\ref{0.4}) are chosen to be these $A_n(\omega).$
Let us define $a_n$ as
\begin{equation}\label{0.2}
a_n:=\int_{-1}^1 f(t) p_n^{(\gamma,\delta)}(t)\rho^{(\gamma,\delta)}(t)dt,\; \gamma,\delta>-1,
\end{equation}
for a suitable function $f$ such that it exists. The $a_n$ are called the Fourier-Jacobi coefficients of the function $f.$
Consider the scalars $d_n$ in the random series (\ref{0.4}) to be  these $a_n.$
The series (\ref{0.4}) is now called the random Fourier-Jacobi series.
It is observed that the different mode of convergence of these random Fourier--Jacobi series depends on the parameters $\gamma,\delta,\eta$ and $\tau$ associated with the Jacobi polynomials, the choice of the scalars $a_n$ and the random variables $A_n(\omega).$
This will create a new area of research in mathematical statistics, opening a scope to look into random
Fourier--Jacobi series of functions from different function spaces, their summability, etc.

This article is organized as follows.
The convergence of the random series (\ref{0.4}) associated with symmetric stable process, for $\gamma,\delta>-1$ is presented in section 3.
In this section, the convergence of the series (\ref{0.4})
for the particular value $\gamma=1/2,\delta=-1/2$ is also established. For this case
$\varphi_n(t)$ are considered to be the weighted orthonormal Jacobi polynomials $u_n^{(1/2,-1/2)}(t),\;t \in [-1,1]$ (see equation (\ref{1.6})).
Further, the random Fourier-Jacobi series (\ref{0.4}) associated with Wiener process is studied in section 4, where
 the orthogonal polynomials $\varphi_n(t)$ are considered to be the modified Jacobi polynomials $q_n^{(\gamma,\delta)}(t)$ (see equation (\ref{1.28})).
All the random Jacobi series dealt with are found to be convergent to the stochastic integral. The continuity property of these sum functions are discussed in section 5.

In the following section we recall some symbols and results on the Fourier-Jacobi series of different classes of functions.

\section{Preliminaries and results on Fourier-Jacobi series}

This section serves to recall the relevant materials from \cite{B}, \cite{PH}, \cite{Z} which will be used in the forthcoming sections.

Zorshchikov \cite{Z} in 1967 first derived the condition which guarantees the uniform convergence of
Fourier-Jacobi series
\begin{equation}\label{1.1}
\mathbf{s}^{(\gamma,\delta)}(f,y):=\sum_{n=0}^\infty  a_n p_n^{(\gamma,\delta)}(y),\; n \in \mathbb{N}
\end{equation}
 on the whole segment $[-1,1],$ where $a_n$ represents the Fourier-Jacobi coefficients
 of continuous functions $f$
defined as in (\ref{0.2}).
He established this result if the scalars $a_n$ are the Fourier-Jacobi coefficients of some class of continuous functions $f$ in $DL[-1,1]$
whose modulus of continuity decreases rapidly in comparison to $\log {n}.$
The class of functions $DL[a,b]$ is defined below.
\begin{definition}[\cite{K}] \label{d2}
The modulus of continuity of a function $f$ in the space of continuous functions on $[a,b]$ with a uniform norm is defined as
\begin{equation*}
\varpi(f,\epsilon,[a,b]):=\sup_{|x-t|\leq \epsilon}\Big\{|f(x)-f(t)|,\;x,t \in[a,b]\Big\}.
\end{equation*}
\end{definition}
\begin{definition}[\cite{K}] \label{d4}
The $DL{[a,b]}$ represents the class of functions $f$ in $C[a,b],$ if
\begin{equation*}
\varpi\Big(f,\frac{1}{n},[a,b]\Big)ln\;n=o(1),\,as \; n\rightarrow \infty.
\end{equation*}
\end{definition}

In 1973, Prasad and Hayashi \cite{PH} considered the weighted Jacobi polynomials $u_n^{(1/2,-1/2)}$ which are
defined by
\begin{equation}\label{1.6}
u_n^{(1/2,-1/2)}(y):={\Bigg[\frac{(2n+1)\Gamma(n+1)\Gamma(n+1)}{2\Gamma(n+3/2)\Gamma(n+1/2)}\Bigg]}^{1/2}p_n^{(1/2,-1/2)}(y),
\end{equation}
where $p_n^{(1/2,-1/2)}$ is the $n$th degree normalized Jacobi polynomial $p_n^{(\gamma,\delta)}$, with $\gamma=1/2$, $\delta=-1/2$,
and investigated the convergence of the Jacobi series
 \begin{equation}\label{1.5}
\sum_{n=0}^\infty b_n u_n^{(1/2,-1/2)}(y).
\end{equation}
They showed that if the scalars $b_n$ defined as
\begin{equation}\label{1.7}
b_n:=\int_{-1}^1 f(t) u_n^{(1/2,-1/2)}(t)\rho^{(1/2,-1/2)}(t)dt
\end{equation}
 are the Fourier-Jacobi
coefficients of some functions $f$ whose $p$th derivative is continuous and belongs to the
Lipschitz class of order $\mu$ less than $1$ on $[-1, 1],$ then the Jacobi series~(\ref{1.5}) with weighted Jacobi polynomials $u_n^{(1/2,-1/2)}(y)$ converges uniformly to $f(y)$ on $[-1,1].$ We will denote the class of function $f$ whose $p$th derivative is continuous and belongs to the Lipschitz class of order $\mu$ less than $1$ on $[-1, 1]$ as $LC^{(p,\mu)}_{[-1,1]}.$

\textit{As we know the Lipschitz class of order $\mu <1$ is the set of functions $f$ satisfy
\begin{equation*}
|f(x)-f(y)| \leq C {|x-y|}^\mu,
\end{equation*}
 for all $x,y \in [-1,1],$ where $C$ is a constant independent of $x,y.$ }

The result for the convergence of the Fourier-Jacobi series of $f$ on the entire segment of orthogonality $[-1,1]$ established by Belen'kii \cite{B} in 1989 is stated below. 
\begin{theorem}\label{t2}
Let $\gamma>-1$ and $\delta> -1.$
If the Fourier-Jacobi series (\ref{1.1}) of a function $f \in DL[-1,1]$ is convergent at $\pm 1,$ then the Fourier-Jacobi series of $f$ converges uniformly on the whole segment of orthogonality $[-1,1].$
\end{theorem}

\section{Random Fourier-Jacobi series associated with symmetric stable process}

Let $X(t,\omega)$ be a symmetric stable process of index $\alpha \in [1,2].$ Consider the random Fourier-Jacobi series
\begin{equation}\label{1.11}
\sum_{n=0}^\infty a_n A_n(\omega) p_n^{(\gamma,\delta)}(y)
\end{equation}
in Jacobi polynomials $p_n^{(\gamma,\delta)}(y), \gamma,\delta>-1$ with random coefficients $A_n(\omega)$ as defined in (\ref{1.13.0}).
Theorem \ref{T1} below establishes the convergence of the series (\ref{1.11}) in probability to the stochastic integral
 \begin{equation}\label{1.14}
\int_{-1}^1 f(y,t)\rho^{(\eta,\tau)}(t)dX(t,\omega),\; \eta,\tau \geq 0,
\end{equation}
depending on the choice of the scalars $a_n$ defined as in (\ref{0.2}).
In fact,
a sequence of random variables $X_n$ is said to converge in probability to a random variable $X$ if
\begin{equation*}
\lim_{n\rightarrow \infty} P(|X_n - X|>\epsilon) = 0, \; \epsilon>0.
\end{equation*}
 The inequality in the following lemma is required to establish this result.
\begin{lemma}[\cite{NPM}] \label{l1}
If $X(t,\omega)$ is a symmetric stable process of index  $\alpha,$ for $1 \leq \alpha \leq 2$ and $f(t)$ is any function in $L^p{[a,b]},\;p\geq 1,$ then for all $\epsilon >0,$
\begin{equation*}
P\Bigg(\Bigg|\int_a^b f(t)dX(t,\omega)\Bigg|> \epsilon\Bigg)
\leq
\frac{C2^{\alpha+1}}{(\alpha+1)\epsilon'^{\alpha}}
\int_a^b|f(t)|^\alpha dt,
\end{equation*}
where $C$ is a positive constant and $\epsilon'<\epsilon.$
\end{lemma}
\begin{theorem}\label{T1}
Let $X(t,\omega),t \in \mathbb{R}$ be a symmetric stable process of index $\alpha \in[1,2],$ and $A_n(\omega)$ be the Fourier-Jacobi coefficients of $X(t,\omega)$ defined as in
(\ref{1.13.0}). If $a_n$ are the Fourier-Jacobi coefficients of some functions $f$ in $DL[-1,1]$ such that
 the Fourier-Jacobi series~(\ref{1.1}) of $f$ is convergent at the end point of segment $[-1,1],$
then the random Fourier-Jacobi series~(\ref{1.11})
converges in probability to the stochastic integral~(\ref{1.14}) on the entire segment $[-1,1],$ for $\eta,\tau \geq 0$ and $\gamma,\delta >-1.$
\end{theorem}
\begin{proof}
Consider the $n$th partial sum of random Fourier-Jacobi series (\ref{1.11})
\begin{align}
\mathbf{S}_n^{(\gamma,\delta)}(f,y,\omega)
&:=\sum_{k=0}^n a_k A_k(\omega)p_k^{(\gamma,\delta)}(y) \label{1.27} \\
&=\sum_{k=0}^n a_k \Bigg(\int_{-1}^1 p_k^{(\gamma,\delta)}(t)\rho^{(\eta,\tau)}(t)dX(t,\omega)\Bigg) p_k^{(\gamma,\delta)}(y) \nonumber \\
&= \int_{-1}^1 \Bigg( \sum_{k=0}^n a_k p_k^{(\gamma,\delta)}(y) p_k^{(\gamma,\delta)}(t)\Bigg)\rho^{(\eta,\tau)}(t)dX(t,\omega) \nonumber \\
&=\int_{-1}^1 \mathbf{s}_n^{(\gamma,\delta)}(f,y,t) \rho^{(\eta,\tau)}(t)dX(t,\omega), \nonumber
\end{align}
where
\begin{equation*}
 \mathbf{s}_n^{(\gamma,\delta)}(f,y,t):=\sum_{k=0}^n a_k p_k^{(\gamma,\delta)}(y) p_k^{(\gamma,\delta)}(t).
\end{equation*}
Since $f \in DL[-1,1]$ are continuous in $[-1,1]$ and $C[-1,1]$ is dense in $L^p{[-1,1]},\;p \geq 1,$
using Lemma \ref{l1},
\begin{align*}
P&\Bigg(\Big|\int_{-1}^1 f(y,t)\rho^{(\eta,\tau)}(t)dX(t,\omega)-\mathbf{S}_n^{(\gamma,\delta)}(f,y,\omega)\Big| >\epsilon\Bigg) \\
&= P\Bigg(\Big|\int_{-1}^1 f(y,t)\rho^{(\eta,\tau)}(t)dX(t,\omega)-\int_{-1}^1 \mathbf{s}_n^{(\gamma,\delta)}(f,y,t)\rho^{(\eta,\tau)}(t)dX(t,\omega)\Big|>\epsilon\Bigg)\\
&\leq \frac{C2^{\alpha+1}}{(\alpha+1)\epsilon'^\alpha}\int_{-1}^1 \Big|\Big(f(y,t)-\mathbf{s}_n^{(\gamma,\delta)}(f,y,t)\Big)\rho^{(\eta,\tau)}(t)\Big|^\alpha dt,
\end{align*}
where $ \epsilon'<\epsilon.$ Since for $f \in DL[-1,1],$ the $n$th partial sum $\mathbf{s}_n^{(\gamma,\delta)}(f,t)$ converges uniformly to $f(t)$ (Theorem \ref{t2}), and
 $\rho^{(\eta,\tau)}(t)$ is bounded in $[-1,1],$ for $\eta,\tau \geq 0,$
 we have
\begin{equation*}
\lim_{n \rightarrow \infty} \int_{-1}^1 \Big|\Big(f(t)-\mathbf{s}_n^{(\gamma,\delta)}(f,t)\Big)\rho^{(\eta,\tau)}(t)\Big|^\alpha dt=0,\; for \; \alpha \in [1,2].
\end{equation*}
This implies the convergence of
the random Fourier-Jacobi series (\ref{1.11})
in probability to the stochastic integral
 (\ref{1.14}).
\end{proof}
The following theorem is concerned about a particular value of $\gamma,\delta$ i.e. $1/2,-1/2$ respectively.
Consider the weighted orthonormal Jacobi polynomial $u_n^{(1/2,-1/2)}(y)$ as defined in (\ref{1.6}) instead of the orthonormal Jacobi polynomial $p_n^{(\gamma,\delta)}(y).$ The random coefficients in the series (\ref{0.4}) are
chosen to be the weighted Fourier-Jacobi coefficients of $X(t,\omega).$ Let us denote it as $B_n(\omega)$ which is defined as
\begin{equation}\label{1.15}
B_n(\omega):=\int_{-1}^1 u_n^{(1/2,-1/2)}(t)\rho^{(\eta,\tau)}(t)dX(t,\omega),\;for \; \eta,\tau \geq 0.
\end{equation}
 The random series that will be dealt with now
 \begin{equation}\label{1.5.1}
\sum_{n=0}^\infty b_nB_n(\omega)u_n^{(1/2,-1/2)}(y).
 \end{equation}
 It is proved that the random Jacobi series (\ref{1.5.1}) in weighted orthonormal Jacobi polynomials
 $u_n^{(1/2,-1/2)}(y)$ converges in probability,
if $b_n$ are the weighted Fourier-Jacobi coefficients of the functions $f$ in the class $LC_{[-1,1]}^{(p,\mu)}$ as defined in (\ref{1.7}).
\begin{theorem}\label{T2}
Let $X(t,\omega),t \in \mathbb{R}$ be a symmetric stable process of index $\alpha \in [1,2]$ and $B_n(\omega)$ be defined as in (\ref{1.15}). If $b_n$ are the Fourier-Jacobi coefficients
of some functions $f$ in the class $LC_{[-1,1]}^{(p,\mu)},$ with $p+ \mu \geq 5/2,$
then the random Fourier-Jacobi series (\ref{1.5.1})
converges in probability to the stochastic integral
(\ref{1.14}), for $\eta,\tau \geq 0.$
\end{theorem}
\begin{proof}
The $n$th partial sum of random Fourier-Jacobi series (\ref{1.5.1}) is
\begin{eqnarray*}
\mathbf{S}_n^{(1/2,-1/2)}(f,y,\omega)&:=&\sum_{k=0}^n b_k B_k(\omega)u_k^{(1/2,-1/2)}(y)\\
&=&\int_{-1}^1 \mathit{v}_n^{(1/2,-1/2)}(f,y,t) \rho^{(\eta,\tau)}(t)dX(t,\omega),
\end{eqnarray*}
where $$\mathit{v}_n^{(1/2,-1/2)}(f,y,t):=\sum_{k=0}^n b_k  u_k^{(1/2,-1/2)}(y) u_k^{(1/2,-1/2)}(t).$$
With the use of inequality in Lemma \ref{l1},
\begin{align}
P&\Bigg(\Big|\int_{-1}^1f(y,t)\rho^{(\eta,\tau)}(t)dX(t,\omega)-\mathbf{S}_n^{(1/2,-1/2)}(f,y,\omega)\Big| >\epsilon\Bigg) \nonumber \\
&\leq \frac{C2^{\alpha+1}}{(\alpha+1)\epsilon'^\alpha}\int_{-1}^1 \Big|\Big(f(y,t)-\mathit{v}_n^{(1/2,-1/2)}(f,y,t)\Big)\rho^{(\eta,\tau)}(t)\Big|^\alpha dt, \label{1.16}
\end{align}
where $0<\epsilon'<\epsilon.$
From Theorem $1$ of \cite{PH},
\[
\Big|f(y)- \mathit{v}_n^{(1/2,-1/2)}(y)\Big| \leq \displaystyle{\frac{C_1^*\ln n}{n^{p+\mu-3/2}}},\; y\in[-1,1], \textup{ for } p+\mu \geq 3/2.
\]
Since $\rho^{(\eta,\tau)}(t)$ is bounded in $[-1,1],$
for $\eta,\tau \geq 0,$ the right hand side of (\ref{1.16}) will tend to $0$ if $p+\mu \geq 5/2$
which completes the proof.
\end{proof}

\section{{Random Fourier-Jacobi series associated with Wiener process}}

Let us consider the Wiener process $W(t,\omega)$ and the $n$th degree polynomial
\begin{equation}\label{1.28}
q_n^{(\gamma,\delta)}(t):=p_n^{(\gamma,\delta)}(2t-1),\; \gamma,\delta>-1,\; n \in \mathbb{N} \cup {0}
\end{equation}
instead of $\varphi_n(t)$ in the random series (\ref{0.4}). These polynomials $q_n^{(\gamma,\delta)}(t)$ are defined in the interval $[0,1]$ and form a complete orthonormal
set in $[0,1]$ with respect to the weight function
\begin{equation*}
\sigma^{(\gamma,\delta)}(t):=(1-t)^\gamma t^\delta,\; \gamma,\delta \geq -1.
\end{equation*}
Call $q_n^{(\gamma,\delta)}(t)$ as the modified Jacobi polynomial.
It is easy to see that Theorem {\ref{t2}} now have the following modified form which will be use in this section.
\begin{theorem}\label{t6}
Let $\gamma,\delta>-1.$ If Fourier-Jacobi series
\begin{equation}\label{1.25}
\sum\limits_{n=0}^\infty c_nq_n^{(\gamma,\delta)}(y)
\end{equation}
 is convergent at $0$ and $1,$ then the Fourier-Jacobi series of $f$ converges uniformly on the whole segment of orthogonality $[0,1],$
for $c_n$ are the Fourier-Jacobi coefficients of some functions $f \in DL[0,1]$ defined as
\begin{equation}\label{1.20.1}
c_n:= \int_0^1 f(t)q_n^{(\gamma,\delta)}(t)\sigma^{(\gamma,\delta)}(t)dt.
\end{equation}
\end{theorem}
It is well known that the stochastic integral
$\int_a^b f(t)dW(t,\omega)$
exists in quadratic mean, for $f \in L^2{[a,b]}$ \cite{NPM}.
As we know, a random sequence $\{X_n\}_{n=0}^\infty$ is said to be convergent to a random variable $X$ in quadratic mean, if
$\lim_{n \to \infty}E \left( |X_n-X|^2 \right)=0.$
This stochastic integral is normally distributed random variable with
mean zero and finite variance $\int_a^b |f(t)|^2 dt,$ if $f(t)$ is a continuous function in $[a,b]$
(c.f. \cite[Lukas,~p.~148]{L}).
The $q_n^{(\gamma,\delta)}(t)\sigma^{(\eta,\tau)}(t)$ remains continuous, for $\eta,\tau \geq 0,$ and hence the stochastic integral
\begin{equation}\label{1.19}
C_n(\omega):=\int\limits_0^1 q_n^{(\gamma,\delta)}(t)\sigma^{(\eta,\tau)}(t)dW(t,\omega),
\end{equation}
with weight $\sigma^{(\eta,\tau)}(t),\eta,\tau \geq 0$ exists in quadratic mean.
The random Jacobi series that we will deal with now is
\begin{equation}\label{1.20}
\sum_{n=0}^\infty c_nC_n(\omega)q_n^{(\gamma,\delta)}(y),\; y\in [0,1],
\end{equation}
where $C_n(\omega)$ are random variables associated with the Wiener process $W(t,\omega)$ defined as in (\ref{1.19}) and the scalars $c_n$ are defined as in (\ref{1.20.1}).
The random coefficients $C_n(\omega)$ are no doubt independent random variables.
We will address
$C_n(\omega)$ and $c_n$ as the modified Jacobi coefficients of $W(t,\omega)$ and function $f$ respectively.
The following theorem proved that the random series (\ref{1.20}) converges in quadratic mean to the stochastic integral
\begin{equation}\label{1.21}
\int_0^1 f(y,t)\sigma^{(\eta,\tau)}(t)dW(t,\omega),
\end{equation}
if $c_n$ are the modified Fourier-Jacobi coefficients of some functions $f$ in $DL[0,1].$
\begin{theorem}\label{T4}
Let $W(t,\omega),t\geq 0$ be the Wiener process
and $C_n(\omega)$ be as defined in (\ref{1.19}).
 If the continuous function $f$ is in $DL[0,1]$ and the Fourier-Jacobi series (\ref{1.25}) is convergent at the end points of the segment $[0,1],$
then the random Fourier-Jacobi series (\ref{1.20}) converges in the entire interval $[0,1]$ to the integral (\ref{1.21})
in quadratic mean, for $\eta,\tau \geq 0$ and $\gamma,\delta >-1.$
\end{theorem}
\begin{proof}
Let
\begin{equation}\label{1.20.2}
\mathbf{T}_n^{(\gamma,\delta)}(f,y,\omega):=
\sum_{k=0}^n c_k
C_k(\omega)q_k^{(\gamma, \delta)}(y),\;\gamma,\delta>-1
\end{equation} be the $n$th partial sum of the random Fourier-Jacobi series (\ref{1.20}).
Putting the expression in (\ref{1.19}) of $C_k(\omega)$ in (\ref{1.20.2}) reduces $\mathbf{T}_n^{(\gamma,\delta)}(f,y,\omega)$ into integral form
\begin{equation*}
\int_0^1 \mathfrak{s} _n^{(\gamma,\delta)}(f,y,t)\sigma^{(\eta,\tau)}(t)dW(t,\omega),\; \eta,\tau \geq 0,
\end{equation*}
where
\begin{equation*}
\mathfrak{s}_n^{(\gamma,\delta)}(f,y,t):=\sum\limits_{k=0}^n c_k q_k^{(\gamma,\delta)}(y)q_k^{(\gamma,\delta)}(t),\; for \; t \in [0,1].
\end{equation*}
For $g \in L^2{[a,b]}$ and $W(t,\omega)$ the Wiener process, we know (Lukacs \cite{L})
\begin{equation}\label{1.22}
E\Big|\int_a^b g(t)dW(t,\omega)\Big|^2=\beta^2\int_a^b |g(t)|^2dt,
\end{equation}
where $\beta$ is a constant associated with the normal law of increment of the process $W(t,\omega),$ for $t \in [a,b].$
Hence
\begin{align}
E&\Big(\Big|\int_0^1 f(y,t)\sigma^{(\eta,\tau)}(t)dW(t,\omega)-T_n^{(\gamma,\delta)}(f,y,\omega)\Big|^2\Big) \nonumber \\
&= E\Big(\Big|\int_0^1 f(y,t)\sigma^{(\eta,\tau)}(t)dW(t,\omega)- \int_0^1 \mathfrak{s}_n^{(\gamma,\delta)}(f,y,t)\sigma^{(\eta,\tau)}(t)dW(t,\omega)\Big|^2\Big) \nonumber \\
&=\beta^2\int_0^1 \Big|\Big(f(y,t)-\mathfrak{s}_n^{(\gamma,\delta)}(f,y,t)\Big)\sigma^{(\eta,\tau)}(t)\Big|^2dt. \label{0.3}
\end{align}
Since $\sigma^{(\eta,\tau)}(t)$ is bounded, for $\eta,\tau \geq 0,$ and $\mathfrak{s}_n^{(\gamma,\delta)}(f,t)$ converges uniformly to $f(t)$ on the whole segment $[0,1],$
if $f\in DL[0,1], \; \gamma,\delta\geq -1$ (Theorem \ref{t6}).
Then the integral~(\ref{0.3}) converges to $0$ implies the convergence of the series~(\ref{1.20}) in quadratic mean to the integral~(\ref{1.21}).
\end{proof}

\section{{Continuity property of the sum functions}}

It is observed that the sum function of the random Fourier-Jacobi series (\ref{0.2}) associated with symmetric stable processes $X(t, \omega)$ of index $\alpha \in [1,2]$ is weakly continuous in probability, where as the sum function of the random series (\ref{1.20}) associated with Wiener process is continuous in quadratic mean.
It is known that a function $f(t,\omega)$ is said to be weakly continuous in probability at $t = t_0,$
if for all $\epsilon > 0,$
$
\lim_{h \rightarrow 0}P(\vert f(t_0 + h, \omega) - f(t_0, \omega) \vert > \epsilon) = 0.
$
If a function $f(t,\omega)$ is weakly continuous at every $t_0 \in[a,b],$ then the function $f(t,\omega)$ is weakly continuous in
probability in the closed interval
$[a, b].$
These facts are established in the following two theorems respectively.
The proofs use the following lemma.
\begin{lemma}[{\cite[p.~37]{ZY}}]  \label{l2}
If $f$ is a periodic function or in $L^p,1 \leq p<\infty$ or a continuous function, then the integral
\begin{equation*}
\Bigg\{\int_a^b \Big|f(x+t)-f(x)\Big|^pdx\Bigg\}^{1/p}
\end{equation*}
tends to $0$ as $t$ tends to $0.$
\end{lemma}
\begin{theorem} \label{T3}
If $X(t,\omega)$ is a symmetric stable process of index $\alpha \in [1,2]$ and $f$ is a function of any class of  $LC_{[-1,1]}^{(p,\mu)}$ or $DL[-1,1],$ with the conditions as stated in Theorem~\ref{T1} and Theorem~\ref{T2} respectively,
then the sum function~(\ref{1.14}) of the random Fourier-Jacobi series~(\ref{1.11}) as well as (\ref{1.5.1}) is weakly continuous in probability, for $\eta,\tau \geq 0.$
\end{theorem}
\begin{proof}
With the help of Lemma \ref{l1},
\begin{equation*}
P\Bigg(\Big|\int_{-1}^1f(x,t)\rho^{(\eta,\tau)}(t)dX(t,\omega)-
\int_{-1}^1f(y,t)\rho^{(\eta,\tau)}(t)dX(t,\omega)\Big|>\epsilon \Bigg)
\end{equation*}
\begin{equation}\label{1.17}
\leq
\frac{C2^{\alpha+1}}{(\alpha+1)\epsilon'^{\alpha}}
\int_{-1}^1\Big|\Big(f(x,t)-f(y,t)\Big) \rho^{(\eta,\tau)}(t)\Big|^\alpha dt,\; \alpha \in [1,2],
\end{equation}
for $0<\epsilon'<\epsilon.$ Since the Jacobi weight $\rho^{(\eta,\tau)}(t)$ is bounded, for $\eta,\tau \geq 0,$
and by Lemma~\ref{l2}, the integral in (\ref{1.17}) converges to $0,$ if $x \rightarrow y.$
This confirms the weakly continuity in probability of the sum function (\ref{1.14}).
\end{proof}
\begin{theorem}
Let $X(t,\omega),t \in \mathbb{R}$ be the Wiener process $W(t,\omega),t\geq 0,$ and $C_n(\omega)$ be the random variables as defined in (\ref{1.19}). If the Fourier-Jacobi series (\ref{1.25}) of the
 continuous function $f$ in $DL[0,1]$ is convergent at the end point of the segment $[0,1],$
then the sum function (\ref{1.21}) is continuous in quadratic mean.
\end{theorem}
\begin{proof}
As for Theorem \ref{t6} the Fourier-Jacobi series (\ref{1.25}) converges in the whole segment $[0,1].$
By Lukacs \cite[p.~147]{L},
\begin{eqnarray*}
&&E\Big(\Big|\int_0^1 f(y,t)\sigma^{(\eta,\tau)}(t)dW(t,\omega)- f(x,t)\sigma^{(\eta,\tau)}(t)dW(t,\omega)\Big|^2\Big)\\
&=&E\Big(\Big|\int_0^1 \Big( f(y,t)-f(x,t)\Big)\sigma^{(\eta,\tau)}(t)dW(t,\omega)\Big|^2\Big)\\
&=& \beta^2 \int_0^1 \Big|\Big(f(y,t)-f(x,t)\Big)\sigma^{(\eta,\tau)}(t)\Big|^2dt.
\end{eqnarray*}
 For $\eta,\tau \geq 0,$
the Jacobi weight $\sigma^{(\eta,\tau)}(t)$ is bounded.
Hence by Lemma \ref{l2},
as $y$ approaches $x,$ the right-hand side tends to zero. This shows that the sum function (\ref{1.21}) is continuous in the quadratic mean, for $\eta,\tau \geq 0$.
\end{proof}

\subsection*{Acknowledgments}
This research work was supported by University Grant Commission (National Fellowship with letter no-F./2015-16/NFO-2015-17-OBC-ORI-33062).

\EditInfo{October 16, 2022}{November 30, 2022}{Serena Dipierro}

\end{paper}